\documentclass[12pt,reqno]{article}

\usepackage[usenames]{color}

\usepackage{latexsym,amssymb,enumerate,amsmath,xcolor}
\usepackage{lscape}
\usepackage{amsthm}
\usepackage{amsmath}
\usepackage{amsfonts}
\usepackage{amssymb}
\usepackage{hyperref} 
\usepackage{amscd}
\usepackage{hyperref} 
\usepackage{mathrsfs}
\usepackage{fancyhdr}
\usepackage{setspace}
\usepackage{color}
\usepackage{enumitem}
\usepackage{longtable}
\usepackage{float}
\usepackage{graphicx}
\usepackage{caption}

\definecolor{webgreen}{rgb}{0,.5,0}
\definecolor{webbrown}{rgb}{.6,0,0}

\usepackage{fullpage}
\usepackage{float}

\usepackage{graphics}
\usepackage{latexsym}
\usepackage{epsf}

\newtheorem{thm}{Theorem}

\newtheorem{rem}{Remark}

\def\d{\displaystyle}

\def \sp{\vspace{ .25cm}}
\def \spz{\vspace{ .5cm}}

\date{}

\begin{document}
\vspace*{1cm}
\begin{center}
{\Large \bf Bell polynomials and generalized \\[.5ex] Laplace transforms} \\
\rm
\vspace{3ex}

{\small by} \; {\large Paolo Emilio Ricci $^{(\dagger)}$} 
\\
\vspace{4ex}
{\footnotesize \em
$^{(\dagger)}$ International Telematic University UniNettuno  \\
Corso Vittorio Emanuele II, 39, 00186 - Roma, Italia \\
}
\spz

\end{center}

\sp
\noindent {\bf Abstract} - \emph {An extension of the Laplace transform obtained by using 
the Laguerre-type exponentials is first shown. Furthermore, the solution of the Blissard 
problem by means of the Bell polynomials, gives the possibility to associate to any numerical 
sequence a Laplace-type transform depending on that sequence.  Computational techniques for 
the corresponding transform of analytic functions, involving Bell polynomials, are derived}.

\vspace{2ex}
{\footnotesize

}

\section{Introduction}\label{sec1}
\setcounter{equation}{0}

It is almost impossible to cite all the contributes of the Laplace transform: 
\begin{eqnarray*}
\begin{array}{c}
{\cal L}(f) := \d \int_0^\infty \exp^{-1}(s\,t) f(t) \, dt = F(s) \,,
\end{array}
\end{eqnarray*}
to the solution of differential equations \cite{gh-os, wid}.
Indeed, the Fourier and the Laplace transforms are the most useful tools in 
Analysis and Mathematical Physics.\\
Actually these transforms are nothing but functions acting in function spaces,
so that it is quite obvious that many others transforms can be defined 
similar to them.

The Bell polynomials \cite{be} have been applied in many different fields of mathematics. 
In order to avoid unuseful repetitions, which would be classified as plagiarism by the 
modern \emph{artificial deficiency}, we limit ourselves to recall the articles 
\cite{be-r, c-r, dc-r, n-r3}.

Some generalized forms of Bell polynomials already appeared 
in literature (see e.g. \cite{fu, r-s}).  The multivariare case was also considered 
in \cite{b-n-r, n-r1, n-r, f-d-d-b}. Connections with number theory have been examined 
in \cite{n-r4, qi}. See also \cite{w-w}.

To the author's knowledge a connection of the two before mentioned topics has not been
considered in literature.  This is the subject of the present article which the proposal 
to introduce a wide extension of the Laplace tranform by using the Blissard umbral calculus.
Computational techniques in case of analytic functions are also given in the last 
sections.

The obtained results, although formal, since they are based on an umbral approach,
allow to consider infinite many other transforms which can be computed by essentially 
algebraic methods.

\section{Recalling the Bell polynomials}\label{sec2} 

Considering the $n$-times differentiable functions $x=g(t)$ and $y=f(x)$, defined in 
given intervals of the real axis, the composite function $\Phi(t):= f(g(t))$, can 
be differentiated with respect to $t$, up to the $n$th order, by using the chain rule.\\
\noindent We use the notations:
$$
\Phi_m:= D_t^m \Phi(t), \qquad f_h:= D_x^hf(x)|_{x=g(t)}, \qquad
g_k:= D_t^k g(t).
$$
Then the $n$th derivative  of $\Phi(t)$ is represented by
$$
\Phi_n = Y_n(f_1,g_1;f_2,g_2;\dots;f_n,g_n),
$$ 
where $Y_n$ denotes the $n$th Bell polynomial.\\
\noindent The first few Bell polynomials are:
\begin{eqnarray}
\label{1}
\begin{array}{l}
Y_1(f_1,g_1)=f_1 g_1 \\
Y_2(f_1,g_1;f_2,g_2)=f_1g_2+f_2g_1^2 \\
Y_3(f_1,g_1;f_2,g_2;f_3,g_3)=f_1g_3+f_2(3g_2g_1)+
f_3g_1^3 \\
\dots
\end{array}
\end{eqnarray}
\noindent Further examples can be found in \cite{rio}, p. 49,
where a recursion formula and the explicit expression given by
the Fa\`a di Bruno formula is also recalled.\\
\noindent A proof of the Fa\`a di Bruno formula based on the 
{\em umbral calculus} is given in \cite{rom} and \cite{r-r}.
However, the Fa\`a di Bruno is not convenient by the computational point 
of view, owing the higher computational complexity with respect to the 
recursion.\\
\noindent The traditional form of the Bell polynomials \cite{co} is given by:
\begin{eqnarray}
\label{2}
\begin{array}{c}
Y_n(f_1,g_1;f_2,g_2;\dots;f_n,g_n) = \d \sum_{k=1}^n B_{n,k}(g_1,g_2,\dots ,g_{n-k+1}) f_k \,,
\end{array}
\end{eqnarray}
where the $B_{n,k}$ satisfy the recursion \cite{co}:
\begin{eqnarray}
\label{3}
      B_{n,k}(g_1,g_2,\dots ,g_{n-k+1}) =
      \d \sum_{h=0}^{n-k} \ {{n-1} \choose {h}} \
      B_{n-h-1,k-1}(g_1,g_2,\dots ,,g_{n-h-k+1})\, g_{h+1} \,.
\end{eqnarray}
The $B_{n,k}$ functions for any $k=1,2,\dots,n$ are polynomials 
in the $g_1,g_2,\dots,g_n$ variables homogeneous of degree $k$ and 
{\em isobaric} of weight $n$ (i.e. they are linear combinations of 
monomials $g_1^{k_1}g_2^{k_2}\cdots g_n^{k_n}$
whose weight is constantly given by $k_1+2k_2+\ldots +nk_n =n$), so that
\begin{eqnarray}
\label{4}
B_{n,k}(\alpha \beta \, g_1, \alpha \beta^2 g_2,\dots , \alpha \beta^{n-k+1} g_{n-k+1})
= \alpha^k \beta^n  B_{n,k}(g_1,g_2,\dots ,g_{n-k+1}) \,,
\end{eqnarray}
and
\begin{eqnarray}
\label{5}
Y_n(f_1, \beta \, g_1;f_2, \beta^2 g_2;\dots;f_n, \beta^n g_n) = \beta^n \, 
Y_n(f_1,g_1;f_2,g_2;\dots;f_n,g_n) \,.
\end{eqnarray}

\section{The Blissard problem}\label{sec3}

John Blissard (1803-1875) published in 1861-1862 papers \cite{bli} introducing a symbolic method
showing that some sequences of numbers $\{b_k\}$ can be substituted by powers $\{b^k\}$ so as to 
obtain valid formulas. The Bernoulli numbers was shown to be a first example of such a sequence. 
The Blissard symbolic method at present is called the \emph{umbral calculus}, a term coined by 
J.J. Sylvester. 

The modern version of the umbral calculus  \cite{r-r, rom} considers the umbral algebra, as the 
algebra of linear functionals on the vector space of polynomials, with the product defined by a 
binomial type formula.  An extensive bibliogaphy of the subject can be found in \cite{dib}.

The so called Blissard problem is described as follows \cite{rio}.\\
Given the formal power series
\begin{eqnarray}
\label{6}
e^{a t} = \sum_{k=0}^\infty \frac{a^k t^k}{k!} =
\sum_{k=0}^\infty \frac{a_k t^k}{k!},
\end{eqnarray}
associted to the sequence $a = \{ a_k \}$, where 
\begin{eqnarray}
\label{7}
a^k := a_k, \qquad \forall k \geq 0, \qquad a_0:= 1,
\end{eqnarray}
the solution of the equation
\begin{eqnarray}
\label{8}
e^{at} e^{bt} = 1
\end{eqnarray}
with respect to the unknown sequence $b = \{ b_n \}$, is given by
\begin{eqnarray}
\label{9}
\left\{
\begin{array}{l}
b_0 := 1 ,  \, \\
\phantom{\rule{1 pt}{26 pt}} 
b_n = Y_n ( -1!, a_1; 2!, a_2; -3!, a_3; \dots; (-1)^n n!, a_n), \quad (\forall \, n>0),
\end{array}
\right.
\end{eqnarray}
where $Y_n$ is the $n$th Bell polynomial \cite{rio}.

\section{A first extension of the Laplace transform}\label{sec4}

In \cite{d-r} the Laguerre-type exponentials has been defined, for every integer $r \geq 1$, 
according to the equation:
\begin{eqnarray}
\label{10}
\begin{array}{c}
e_r(x):= \d \sum_{k=0}^\infty \frac{x^k}{(k!)^{r+1}} \, .
\end{array}
\end{eqnarray}
Obviously, for $r=0$, it results: $e_0(x) \equiv \exp(x)$.\\
Actually these functions are a particular case of the Le Roy functions \cite{ler}, 
and more generally of the generalized Mittag–Leffler functions defined in \cite{g-p}, 
and deeply studied in \cite{g-r-m}. 

A comparison among the functions $[e_r(t)]^{-1}$, $(r=1,2)$, 
$\exp(-t)$ and the limit value $\d \lim_{r \rightarrow + \infty} [e_r(t)]^{-1} = 1/(1+t)$ 
shows that for every $r\geq 1$ the functions $[e_r(t)]^{-1} \in L^1(0, +\infty)$, while
the limit value $1/(1+t)$ does not satisfy this condition.
 
Consider the following transforms:
\begin{eqnarray}
\label{11}
\begin{array}{c}
{\cal L}_{1}(f) := \d \int_0^\infty [e_1(s\,t)]^{-1} f(t) \, dt = \d \int_0^\infty 
\left[ \d \sum_{k=0}^\infty \frac{{(s\,t)}^k}{(k!)^2} \right]^{-1} f(t) \, dt = F_1(s) \,,
\end{array}
\end{eqnarray}

\noindent and in general:
\begin{eqnarray}
\label{12}
\begin{array}{c}
{\cal L}_{r}(f) := \d \int_0^\infty [e_r(s\,t)]^{-1} f(t) \, dt = \d \int_0^\infty 
\left[ \d \sum_{k=0}^\infty \frac{{(s\,t)}^k}{(k!)^r} \right]^{-1} f(t) \, dt  = F_r(s) \, .
\end{array}
\end{eqnarray}

As in the ordinary Laplace transform, the integrals in equations (\ref{11})-(\ref{12}) exist
for all real numbers $\d \mathrm{Re}(s)>a$, where the constant $a$, called the convergence 
abscissa,  depends on the function $f$ and determines the region of convergence. \\
Note that the increasing behaviour of the Laguerre-type exponentials, in the interval $(0, +\infty)$ 
is lower with respect to the ordinary exponential, so that for any fixed $f$, we can choose, 
at least, the same convergence abscissa of the ordinary Laplace transform.

For every $r \geq 1$, an approximation of the transform (\ref{12}) is obtained by using 
the truncated Laguerre-type exponential of order $r$, putting for a fixed integer $n$:
\begin{eqnarray}
\label{13}
\begin{array}{c}
\d \int_0^\infty \left[ \d \sum_{k=0}^n \frac{{(s\,t)}^k}{(k!)^r} \right]^{-1} f(t) \, dt 
= F_r^{[n]}(s) \,.
\end{array}
\end{eqnarray}

Since for $r=0$ the Laguerre-type exponentials give back the ordinary 
exponential function, then, in this case, equation (\ref{11}) reduces to the ordinary Laplace
transform, that is it results:
\begin{eqnarray}
\label{14}
\begin{array}{c}
{\cal L}_{0}(f) = {\cal L}(f) := \d \int_0^\infty \exp^{-1}(s\,t) f(t) \, dt = F_0(s) \, .
\end{array}
\end{eqnarray}

\subsection{The inversion formula}

The Laguerre-type exponentials are monotonic increasing functions in the interval 
$(0, +\infty)$, so that can be inverted in the same interval.  Therefore, we can \emph{conjecture} 
that the transform (\ref{12}) (in particular (\ref{11})), admits the inversion formula:
\begin{eqnarray}
\label{15}
\begin{array}{c}
{\cal L}^{-1}_{r}(f)(t) := \d \frac{1}{2 \pi {\rm i}} \d \lim_{\tau \rightarrow \infty} 
\int_{\gamma - {\rm i} \, \tau}^{\gamma + {\rm i} \, \tau} [e_r(s\,t)] F_r(s) \, ds \,,
\end{array}
\end{eqnarray}
where $\gamma$ is a real number so that the contour path of integration is in the 
region of convergence of $F_r(s)$. It should also be possible to transform the contour 
into a closed curve, allowing the use of the residue theorem.

However, the proof of the equation (\ref {15}) is not easy to carry out, since it would be 
necessary to introduce an extension of the Fourier transform based on Laguerre exponentials, 
a topic still far from being obtained.

\sp
\section{The isomorphism \ ${\cal T}_s$ and its iterations}\label{sec5}

In previous articles (see e.g. \cite{r-t}), it was shown that there exist a differential
isomorphism  ${\cal T}:= {\cal T}_s$, acting into the space 
${\cal A}:= {\cal A}_s$ of analytic functions of the 
variable $s$, by means of the correspondence:
\begin{eqnarray*}
D_s:= D \equiv \frac{d}{ds} \; \rightarrow \; {\hat D}_L:= D_s s D_s ; \qquad  s
\cdot \; \rightarrow \; {\hat D}_s^{-1}, \qquad
\end{eqnarray*}
where
\begin{eqnarray*}
{\hat D}_s^{-n}F(s):= \frac{1}{(n-1)!}\int_0^s (s -
\xi)^{n-1}F(\xi) d\xi \,.
\end{eqnarray*}
The isomorphism  ${\cal T}:= {\cal T}_s$  can be iterated
producing a set of generalized Laguerre derivatives as follows.
According to the results in \cite{r-t} we put, for every integer $m \geq 1$,
\begin{eqnarray*}
{\cal T}_s^{m-1}{\hat D}_L = {\cal T}_s^{m-1} (DsD) = DsDsD \cdots
sD =: {\hat D}_{mL},
\end{eqnarray*}
where the last operator contains  $s+1$  ordinary derivatives,
denoted by $D \equiv D_s$. 

The action of ${\cal T}_s$, on powers, and
consequently on all functions belonging to \ ${\cal A}:= {\cal
A}_s$ is as follows:
\begin{eqnarray*}
{\hat D}_s^{-n} (1)= \frac{s^n}{n!} \ ,
\end{eqnarray*}
and, by induction:
\begin{eqnarray*}
{\cal T}_s^{m-1} {\hat D}_s^{-1}(1) = {\hat D}_{{\cal
T}_s^{m-1}}^{-1} (1)  \quad \Rightarrow \quad {\hat D}_{{\cal
T}_s^{m-1}}^{-n} (1) = \frac{s^n}{(n!)^s} \ . \quad
\nonumber 
\end{eqnarray*}

Note that the Laguerre-type exponentials are obtained, acting with these 
iterated isomorphisms on the classical exponential, since:
\begin{eqnarray*}
{\cal T}_s^m (e^s) = \sum_{k=0}^\infty \frac{{\cal T}_s (s^k)
}{(k!)^m} = \sum_{k=0}^\infty \frac{s^k}{(k!)^{m+1}} = e_m(s) \,.
\end{eqnarray*}

It has been shown in a number of articles \cite{b-d-r, be-br-ri, br-ce-ri, ce-ge-ri}, 
that new sets of special functions, namely the Laguerre-type special functions, can be 
introduced and some of their applications have been considered in \cite{br-ri, da-he-ri, 
dea-ri, m-r3}.

\subsection{Computation via the isomorphisms ${\cal T}_s$}

Acting with the isomorphism ${\cal T}_s$ on both sides of equation (\ref{14}), we find:
\begin{eqnarray*}
\begin{array}{c}
\d \int_0^\infty {\cal T}_s [\exp^{-1}(s\,t)]\, f(t) \, dt = {\cal T}_s [F_0(s)] \,,
\end{array}
\end{eqnarray*}
that is
\begin{eqnarray*}
\begin{array}{c}
\d \int_0^\infty [e_1(s\,t)]^{-1} f(t) \, dt  = {\cal T}_s [F_0(s)] \,,
\end{array}
\end{eqnarray*}
so that, comparing this result with equation (\ref{11}), we find:
\begin{eqnarray*}
\begin{array}{c}
{\cal T}_s [F_0(s)] = F_1(s) \,.
\end{array}
\end{eqnarray*}

Of course this equation can be generalized starting from (\ref{12}), obtaining:
\begin{eqnarray*}
\begin{array}{c}
\d \int_0^\infty {\cal T}_s [e_r(s\,t)]^{-1} f(t) \, dt = \d \int_0^\infty 
\d \sum_{k=0}^\infty [e_{r+1}(s\,t)]^{-1} f(t) \, dt  = {\cal T}_s F_r(s) \,,
\end{array}
\end{eqnarray*}
and therefore
\begin{eqnarray*}
\begin{array}{c}
{\cal T}_s [F_r(s)] = F_{r+1}(s) \,.
\end{array}
\end{eqnarray*}

\section{A more general extension of the Laplace transform}\label{sec6}

A further extension of the transforms (\ref{11})-(\ref{12}) is as follows.\\
Given the sequence $a := \{ a_k \} = (1, a_1, a_2, a_3, \dots )$, we consider the function:
\begin{eqnarray}
\label{16}
\begin{array}{c}
\d \frac{1}{1 + a_1 t + a_2 \frac{t^2}{2!} + a_3 \frac{t^3}{3!} + \dots } \, \quad (t\geq 0).
\end{array}
\end{eqnarray}
When $a_k = 1/(k!)^{r}$ the function (\ref{10}) is recovered, and for $r=0$ we find 
again $\exp(-t)$.

Note that the functions (\ref{16}) are complete monotonic functions decreasing from
the initial value $1$, at $t=0$, and vanishing at infinity.

Therefore, according to the umbral method, we put by definition:
\begin{eqnarray}
\label{17}
\begin{array}{c}
{\cal L}_a(f):= \d \int_0^\infty \frac{f(t)}{\d \sum_{k=0}^\infty \frac{a_k (s\,t)^k}{k!}}\, dt = 
\d \int_0^\infty \frac{f(t)}{\d \sum_{k=0}^\infty \frac{a^k (s\,t)^k}{k!}}\, dt =
F_a(s)
\end{array}
\end{eqnarray}

Recalling the Blissard problem, the solution of the the umbral equation
\begin{eqnarray}
\label{18}
\begin{array}{c}
\d \frac{1}{\d \sum_{k=0}^\infty \frac{a^k (s\,t)^k}{k!}} = 
\d \sum_{k=0}^\infty \frac{b^k (s\,t)^k}{k!}\, 
\end{array}
\end{eqnarray}
that is
\begin{eqnarray}
\label{19}
\begin{array}{c}
\exp[a (s\,t)] \exp[b (s\,t)] = 1 \,,
\end{array}
\end{eqnarray}
is given by equation (\ref{9}).

Therefore, the generalized Laplace transform (\ref{13}) writes:
\begin{eqnarray}
\label{20}
\begin{array}{c}
{\cal L}_a(f) := \d \int_0^\infty f(t) \left[ 1+ \d \sum_{k=1}^\infty 
Y_k (-1!, a_1; 2!, a_2; \dots ; (-1)^k k!, a_k) \frac{(s\,t)^k}{k!} \right] \, dt = F_a(s) \,.
\end{array}
\end{eqnarray}
By using equation (\ref{2}), equation (\ref{20}) becomes
\begin{eqnarray}
\label{21}
\begin{array}{c}
{\cal L}_a(f) = \d \int_0^\infty f(t) \left[ 1+ \d \sum_{k=1}^\infty \sum_{h=1}^k (-1)^h h! \,
B_{k,h}(a_1, a_2, \dots, a_{k-h+1}) \d \frac{(s\,t)^k}{k!} \right] \, dt = \\
\phantom{\rule{1 pt}{30 pt}}
= \d \int_0^\infty f(t) \left[ 1+ \d \sum_{k=1}^\infty \sum_{h=1}^k (-1)^h h! \,
B_{k,h}(a_1, a_2, \dots, a_{k-h+1}) \, \d \frac{(s\,t)^k}{k!}  \right] \, dt = F_a(s) \,.
\end{array}
\end{eqnarray}

It is convenient to introduce the definition
\begin{eqnarray}
\label{22}
\begin{array}{c}
C_k(a):= \d \sum_{h=1}^k (-1)^h h! \, B_{k,h}(a_1, a_2, \dots, a_{k-h+1}) \,, \quad C_0(a):= 1 \,,
\end{array}
\end{eqnarray}
so that equation (\ref{20}) writes
\begin{eqnarray}
\label{23}
\begin{array}{c}
{\cal L}_a(f) = \d \int_0^\infty f(t) \d \sum_{k=0}^\infty C_k(a) \, \d \frac{(s\,t)^k}{k!} \, dt = F_a(s) \,.
\end{array}
\end{eqnarray}

Note that the ordinary Laplace transform corresponds to the sequence:
$a = (1,1,\dots, 1, \dots)$, that is $a_k \equiv 1\,, \forall k\geq 1$.  Therefore, it results:
\begin{eqnarray}
\label{24}
{\cal L}_{(1,1,\dots, 1, \dots)}(f) \equiv {\cal L}(f) \,.
\end{eqnarray}
In this case, we find:
\begin{eqnarray*}
B_{k,h} (1,1,\dots ,1) = S(k,h) \,, 
\end{eqnarray*}
where $S(k,h)$ are the Stirling numbers of the second kind \cite{co}.\\
Then
\begin{eqnarray*}
C_k (1,1,\dots ,1) = \d \sum_{h=1}^k (-1)^h h! \, S(k,h)  \,, \quad C_0(a):= 1 \,.
\end{eqnarray*}
Recalling the known identity
\begin{eqnarray*}
\d \sum_{h=1}^k (-1)^{k-h} h! \, S(k,h) = 1 \,, 
\end{eqnarray*}
we find
\begin{eqnarray*}
C_k (1,1,\dots ,1) = (-1)^k \,, \quad \forall k \geq 0 \,,
\end{eqnarray*}
so that the ordinary expression of the Laplace transform is recovered.

\sp
A number of sums defining the coefficients $C_k(a)$ corresponding to different sequences 
$a=\{a_k \}$ can be found in \cite{qi}, however, in what follows, we will assume the 
fundamental hypothesis that preservs the property of the ordinary Laplace transform:\\

\textbf{HP.} \emph{For every fixed} $s$ \emph{in the region of convergence, the power series }
$\sum_{k=0}^\infty C_k(a) \, (s\,t)^k/k!$, \emph{in equation} (\ref{23}) \emph{has an exponential 
decay to zero when} $t \rightarrow \infty$.\\

In this framework, another possibility is to assume $a_k = k!$. In this case equation (\ref{16}) becomes:
\begin{eqnarray}
\label{25}
\begin{array}{c}
\d \frac{1}{1 + t + t^2 + t^3 + \dots } \, \quad (t\geq 0).
\end{array}
\end{eqnarray}

The truncation of the geometric series at the denominator in equation (\ref{25}) produces  
graphs corresponding to the sequences 
$$(1,1,0,0,0,\dots), \ (1,1,1,0,0,0,\dots), \ (1,1,1,1,0,0,0,\dots).$$ 
The decreasing character of the corresponding graphs increases as the number of units  increase.

\begin{rem} {\rm
Increasing the values of the sequence $\{a_k \}$ in equation (\ref{16}), the corresponding
graphs exibhit a more fast decreasing character.
Then the transforms coresponding to the relevant truncations can be limiteded to a small 
interval of the type $[0,L]$, as the values of the function (\ref{16}) become 
negligible outside this interval.}
\end{rem}

\section{Properties}\label{sec7}

From the definition (\ref{23}) the following properties are derived:\\
$\bullet$ \textbf{Linearity}
\begin{eqnarray}
\label{26}
\begin{array}{c}
{\cal L}_a (A f_1 + B f_2 ) = A\,{\cal L}_a (f_1) + B\,{\cal L}_a (f_2) \,.
\end{array}
\end{eqnarray}
$\bullet$ \textbf{Homothetic property}\\
Putting: $x\,a :=(x a_1, x^2 a_2, \dots, x^n a_n , \dots)$, by using the isobaric property (\ref{5})
of the Bell polynomials, it results:
\begin{eqnarray}
\label{27}
\begin{array}{c}
{\cal L}_{x\,a}(f) = F_a(x\,s) \,.
\end{array}
\end{eqnarray}
This can be interpreted as an homothety between the space of the $a$ parametes and that of the 
variable $s$.

$\bullet$ \textbf{Scaling property}
\begin{eqnarray}
\label{28}
\begin{array}{c}
{\cal L}_{a}(f(d\, t)) = \d \frac{1}{d} F_a\left(\frac{s}{d}\right) \,.
\end{array}
\end{eqnarray}

\textbf{Proof.} From equation (\ref{23}) we find:
\begin{eqnarray*}
\begin{array}{c}
F_a \left(\frac{s}{b}\right) = \d \int_0^\infty f(t) \d \sum_{k=0}^\infty C_k(a) \, 
\left(\frac{s}{b}\right)^k \frac{t^k}{k!} \, dt \,,
\end{array}
\end{eqnarray*}
changing variable, putting $t = b\,x$, it results
\begin{eqnarray*}
\begin{array}{c}
F_a \left(\frac{s}{b}\right) = b \,\d \int_0^\infty f(b\, x) \d \sum_{k=0}^\infty C_k(a) \, 
\frac{(s\, x)^k}{k!}\, dx = {\cal L}_{a}(f(b\, x)) \,,
\end{array}
\end{eqnarray*}
that is equation (\ref{28}), up to the change of name of the variable $t$.

$\bullet$ \textbf{Action on the derivative}
\begin{eqnarray}
\label{29}
\begin{array}{c}
{\cal L}_a(f') = \d \int_0^\infty f'(t) \d \sum_{k=0}^\infty C_k(a) \, \frac{(s\, t)^k}{k!} \, dt = \\
\phantom{\rule{1 pt}{26 pt}}
= -\, s \d \int_0^\infty f(t) \d \sum_{k=0}^\infty C_{k+1}(a) \, \frac{(s\, t)^k}{k!} \, dt
- f(0) \,.
\end{array}
\end{eqnarray}
\textbf{Proof.} It is sufficient to integrate by parts and to use the above HP.

\section{Computational techniques}\label{sec8}

According to the above definitions, it is possible to prove the theorems
\begin{thm}
Let $f(t)$ be an analytic function on the real axis. 
Using the Taylor expansion of the function $f(t)$, centered at the origin:
\begin{eqnarray}
\label{30}
\begin{array}{c}
f(t) = \d \sum_{k=0}^\infty c_k \frac{t^k}{k!} \,,
\end{array}
\end{eqnarray}
equation \emph{(\ref{23})} writes:
\begin{eqnarray}
\label{31}
\begin{array}{c}
{\cal L}_a(f) = \d \int_0^\infty \d \sum_{n=0}^\infty \d \sum_{k=0}^n {n \choose k} 
c_{n-k} \, C_k(a) \, \frac{(s\,t)^k}{k!} \, dt = F_a(s) \,.
\end{array}
\end{eqnarray}
\end{thm}
\textbf{Proof. }- In fact, from equations (\ref{23})-(\ref{30}), by using
the Cauchy product of power series, we find
\begin{eqnarray*}
{\cal L}_a(f) = \d \int_0^\infty \sum_{k=0}^\infty c_k \frac{t^k}{k!} 
\d \sum_{k=0}^\infty \, C_k(a) \frac{(st)^k}{k!} \, dt =  \d \int_0^\infty 
\d \sum_{n=0}^\infty \d \sum_{k=0}^n {n \choose k} 
c_{n-k}\, C_k(a) \frac{(s\,t)^k}{k!}\, dt \,,
\end{eqnarray*}
that is the result.

\begin{thm}
Let $f(t)$ be a function expressed by the Laurent expansion:
\begin{eqnarray}
\label{32}
\begin{array}{c}
f(t) = \d \sum_{k=0}^\infty c_k \frac{t^{-k}}{k!} \,,
\end{array}
\end{eqnarray}
then, equation \emph{(\ref{23})} writes:
\begin{eqnarray}
\label{33}
\begin{array}{c}
{\cal L}_a(f) = \d \int_0^\infty \d \sum_{n=0}^\infty \d \sum_{k=0}^n {n \choose k} 
c_{n-k} \, C_k(a) \,  \frac{ s^k \, t^{-n+2k}}{n!}\, dt 
= F_a(s) \,.
\end{array}
\end{eqnarray}
\end{thm}
\textbf{Proof. }- In fact, considering the Cauchy product:
\begin{eqnarray}
\label{34}
\begin{array}{c}
\d \sum_{k=0}^\infty c_k x^{-k} \sum_{k=0}^\infty a_k x^k = \sum_{n=0}^\infty 
\sum_{k=0}^n c_{n-k} x^{-n+k} a_k x^k = \sum_{n=0}^\infty \sum_{k=0}^n a_k c_{n-k} x^{-n+2k} \,.
\end{array}
\end{eqnarray}
from equations (\ref{23})-(\ref{32}), we find
\begin{eqnarray*}
\begin{array}{c}
{\cal L}_a(f) = \d \int_0^\infty \sum_{k=0}^\infty c_k \frac{t^{-k}}{k!} 
\d \sum_{k=0}^\infty C_k(a) \,  \frac{(s\,t)^k}{k!} \, dt = 
\d \int_0^\infty \d \sum_{n=0}^\infty \d \sum_{k=0}^n  \frac{c_{n-k}}{(n-k)!} \frac{C_k(a)}{k!} 
\, s^k t^{-n+2k} \, dt \,,
\end{array}
\end{eqnarray*}
that is the result.

\section{A general isomorphism \ ${\cal T}_s(a)$}\label{sec9}

The results of Section 5 suggests the possibility to introduce a more
general isomorphisms.  The isomorphism ${\cal T}_s$, defined in Section 5 
is determined by the sequence $a^{[1]}:= (1, 1/2!, 1/3!, \dots )\,$. 
Now, given a sequence of nonvanishing real numbers 
$a:= (a_1, a_2, a_3, \dots )$, $(a_k \neq 0, \forall k)$, we can define a correspondence 
acting into the space ${\cal A}:= {\cal A}_s$ of analytic functions of the $s$ variable 
by means of the position:
\begin{equation*}
{\cal T}_s(a) s^n = a_n s^n \,.
\end{equation*}
In particular, the isomorphism ${\cal T}_s$ is recovered, since it results: 
${\cal T}_s(a^{[1]}) \equiv {\cal T}_s$.\\
Even if this isomorphism is not derived form a differential operator, it is still
possible to apply it to the generalized Laplace transform (\ref{23}), obtaining the 
equation:
\begin{equation*}
{\cal T}_s(a) [F_a(s)] = \d \int_0^\infty f(t) \d \sum_{k=0}^\infty C_k(a) 
\, a_k \, (s\,t)^k \, dt \,.
\end{equation*}

\section{Conclusion}\label{sec10} 

\noindent It has been shown that, by exploiting Laguerre-type exponentials, it is possible 
to introduce generalized forms of the Laplace transform that it is supposed to be applied 
in the treatment of differential equations that use the Laguerre derivative instead of the 
ordinary one. For these transforms it was also possible to deduce the transformed functions 
by means of a differential isomorphism studied in previous articles.

The particular form of Laguerre exponentials also suggested a wider extension of the Laplace 
transform, which is associated with a sequence of numbers denoted by the umbral symbol $a$. 
This extension, through the solution of the Blissard problem, which uses Bell's polynomials 
in a natural way, has made it possible to formally define a whole class of transforms, each 
of which is associated with a fixed sequence.  Some fundamental calculation rules have been 
demonstrated for all the generalized transformations considered.

Numerous problems remain open, first of all the existence and the analytical proof of the 
inverse transformation, which should be based on the extension of the Fourier transform to 
the Laguerrian case, with all the problems related to the study of a completely new Fourier 
type analysis.\\
\sp
\noindent \textbf{Acknowledgment} - I want to thank Prof. Dr. Francesco Mainardi for let me know 
his work on the generalized Mittag-Leffler functions and their applications in the Laplace 
transform theory.

\spz
\noindent \textbf{\large Compliance with ethical standards}

\noindent \textbf{Funding} \ {This research received no external funding.}

\noindent \textbf{Conflict of interest} \ The author declares that he has 
not received funds from any institution.

\sp

\end{document}